\documentclass{amsart}
\usepackage{graphicx}
\usepackage{amssymb}
\usepackage{amsmath}
\usepackage{setspace}
\newtheorem{thm}{Theorem}[]
\newtheorem{cor}[]{Corollary}
\newtheorem{lem}[]{Lemma}
\newtheorem{prop}[]{Proposition}
\theoremstyle{definition}

\theoremstyle{remark}
\newtheorem{rem}[]{Remark}
\theoremstyle{definition}

\newcommand{\tr}{{\rm tr}}

\newcommand{\la}{\lambda}

\newcommand{\calD}{\mathcal{D}}

\newcommand{\frakm}{\mathfrak{m}}

\newcommand{\bB}{\mathbb{B}}

\newcommand{\bS}{\mathbb{S}}

\newcommand{\mX}{\mathbf X}

\newcommand{\mU}{\mathbf U}

\newcommand{\mA}{\mathbf A}
\newcommand{\mB}{\mathbf B}
\newcommand{\mC}{\mathbf C}

\newcommand{\mI}{\mathbf I}
\newcommand{\mD}{\mathbf D}

\newcommand{\MW}{\mathbf{W}}

\newcommand{\eqdist}{\ensuremath{\stackrel{\mbox{\upshape\tiny $\calD$}}{=}}}

\author{Tamer Oraby}
\thanks{Department of Mathematical Sciences,
University of Cincinnati, 2855 Campus Way, PO Box 210025,
Cincinnati, OH 45221-0025, USA. E-mail: orabyt@math.uc.edu
Tel:513-5568471}

\title[Limiting Spectra of Girko's Block-Matrix]{The Limiting Spectra of Girko's Block-Matrix}
\begin{document}
\maketitle
\begin{abstract}
To analyze the limiting spectral distribution of some random
block-matrices, Girko \cite{Girko2000} uses a system of canonical
equations from \cite{Girko98}. In this paper, we use the method of
moments to give  an integral form for the almost sure limiting
spectral distribution of such matrices.
\end{abstract}

\section{Introduction and main result}

A random block-matrix is a matrix whose entries are random matrices.
In \cite{Girko98}, Girko studied the spectra of large dimensional
random block-matrices by introducing a system of equations, called
the system of canonical equations, to analyze the spectra. This
system of canonical equations was used later by Girko
\cite{Girko2000} to study a model for which the system is solvable.
The model studied there has many restrictive conditions.

In the current paper, we are going to study the same model under
different conditions for the blocks. The main tool of the proof is
the method of moments. We will follow the proof of the main theorem
by propositions, as applications to the theorem, in which the blocks
are made of some known ensembles like the Gaussian unitary ensemble
and the Wishart random matrix. Free probability theory is used to
prove these propositions.

The spectral measure of an $n \times n$ Hermitian matrix $\mA$ is
$$\mu_{\mA}=\frac1n \sum_{i=1}^n \delta_{\lambda_i},$$ where
$\la_1\leq\la_2\leq\cdots\leq\la_n$ are the eigenvalues of $\mA$. In
this paper we consider random matrices, \emph{i.e.}, matrices where
the entries are random variables on some probability space. In this
case, $\mu_{\mA}$ is, of course, a random measure.

We will denote the weak convergence of probability measures by
$\lim_{n \to \infty} \mu_n \eqdist \mu$ or
\begin{equation*}\label{dist1}
\mu_n \xrightarrow{\calD} \mu \: \mbox{ as $n\to \infty$} .
\end{equation*}

If the moments of measures converge,
\begin{equation*}\label{moment2}
\lim_{n \to \infty} \int_{\mathbb{R}}x^k\mu_n(dx)=
\int_{\mathbb{R}}x^k\mu(dx)
\end{equation*}
for all $k\geq 1$, we will write
\begin{equation*}\label{moment1}
\mu_n \xrightarrow{\frakm} \mu \: \mbox{ as $n\to \infty$}.
\end{equation*}
If $\{\mu_n\}$ is a sequence of random measures which converges in
one of the above senses
 almost surely, we will append the abbreviation "\emph{a.s.}"
 to the above notation. We note that for
$k\geq 1$, the $k^{th}$ moment of $\mu_{\mA}$ is
\begin{equation}\label{matmoment}
\int_{\mathbb{R}} x^k \mu_{\mA}(dx)=\tr_n(\mA^k),
\end{equation}
where $\tr_n(\mA):=\frac1n \sum_{i=1}^n A_{ii}$.

The Kronecker product $\otimes$ of two matrices
$\mA=(a_{ij})_{i,j=1}^k$ and $\mB=(b_{ij})_{i,j=1}^n$ is defined to
be the $nk \times nk$ matrix given by
$$\mA \otimes \mB=
\left[ {\begin{array}{*{20}c}
   a_{11}\mB & a_{12}\mB & {\dots} & a_{1k}\mB  \\
   a_{21}\mB & a_{22}\mB & {\dots} & a_{2k}\mB  \\
   {\vdots} &  {\vdots} & {\ddots} &{\vdots}   \\
   a_{k1}\mB & a_{k2}\mB &  {\dots} & a_{kk}\mB  \\
 \end{array} } \right].
$$
Among the properties of the Kronecker product, we will need the
identity $$\tr_{kn}(\mA \otimes \mB)=\tr_k(\mA) \tr_n(\mB).$$ Also,
if $\mA$ and $\mC$ are two $k \times k$ matrices and $\mB$ and $\mD$
are two $n \times n$ matrices, then $(\mA \otimes \mB)(\mC \otimes
\mD)=\mA \mC \otimes \mB \mD$. Finally, $\mI_k$ is the $k \times k$
identity matrix.

Now we are ready to state the main theorem.
\begin{thm}\label{main}
For $n,k\geq 1$, consider the double array of random block-matrices
$\{\bB_{n,k}\}$ whose terms are given by $\bB_{n,k}=\mI_k \otimes
\mA_n + \MW_k \otimes \mB_n$, where for $n\geq 1$, the matrices
$\mA_n$, $\mB_n$, and $\MW_n$ are Hermitian random matrices of order
$n$, and satisfy the following hypotheses:
\begin{enumerate}
\item \label{1)} There exists a compactly supported
probability measure $\mu_{\omega}$ such that
\begin{equation}\label{main1}
\mu_{\MW_n}\xrightarrow{\frakm} \mu_{\omega} \: \mbox{ as $n\to
\infty$} \qquad a.s.
\end{equation}
\item \label{2)}
For real $t$, there exist probability measures $\psi(t;.)$ such that
\begin{equation*}\label{main2}
\mu_{\mA_n+t\mB_n}\xrightarrow{\frakm} \psi(t;.) \: \mbox{ as $n\to
\infty$} \qquad a.s.
\end{equation*}
and $\psi(t;.)$ has a support that is uniformly bounded for $t$ in
any compact subset of $\mathbb{R}$. \tolerance=1000
\end{enumerate}
\tolerance=500

Under these conditions we have
\begin{equation}\label{main3}
\lim_{n \to \infty} \lim_{k \to \infty}\mu_{\bB_{n,k}} \eqdist
\lim_{k \to \infty} \lim_{n \to \infty}\mu_{\bB_{n,k}} \eqdist \nu
\qquad a.s.,
\end{equation}
where the probability measure $\nu$ is defined as
\begin{equation}\label{main4}
\nu(dx)=\int_{\mathbb{R}} \psi(t;dx) \: \mu_{\omega}(dt).
\end{equation}
\end{thm}
\begin{rem}\label{girko}\noindent
Since the support of $\psi(t;dx)$ is uniformly bounded for $t$ in
$\mbox{supp}(\mu_\omega)$ (the support of $\mu_{\omega}$), then the
probability measure $\nu(dx)$, introduced in \eqref{main4}, is
compactly supported.
\end{rem}

As mentioned above, matrices  of the form $\mI_k \otimes \mA_n +
\MW_k \otimes \mB_n$ were analyzed by Girko \cite{Girko2000}. In
\cite[Theorem 3]{Girko2000}, Girko assumes that $\{\mA_n\}$ and
$\{\mB_n\}$ are two sequences of real symmetric non-random matrices,
$\mB_n$ is a positive definite matrix for each $n \geq 1$, and the
entries of the symmetric matrix $\MW_k$  are independent $\pm1$ with
probability $\frac12$. He  shows under these assumptions that the
spectral probability distribution
$\widehat{F}_{\bB_{n,k}}(x):=\mu_{\bB_{n,k}}\left((-\infty,x]\right)$
of the sequence of random block-matrices $\{\bB_{n,k}\}$ converges,
for almost all $x$'s and with probability one, as both $k$ and $n$
go to infinity to a non-random distribution function that follows
from a complicated equation given in \cite{Girko2000}.

In Theorem \ref{main}, our assumptions allow us to identify the
limit. In the course of our proof, we are also able to derive
Girko's SS-Law (a sum of semi-circular law), see Proposition
\ref{commute}.

\section{Proof of Theorem {\protect\ref{main}}}\label{Proof}

We need the following lemma.
\begin{lem}\cite[p.94, Lemma 3.1]{ferreira}\label{poly}
Fix $k\in \mathbb{N}$, let $T=\{t_0,t_1,\ldots,t_k\}$ be a set of
distinct points in $\mathbb{R}$ and
$P_{n}(t)=a_{0,n}+a_{1,n}t+\cdots+a_{k,n}t^k$ be a polynomial with
$a_{i,n} \in \mathbb{C}$ for every $i$ and $n$. If $P_{n}(t)$
converges for every $t\in T$ as $n \to \infty$, then the limit is a
polynomial of degree $\leq k$, say it is
$P(t)=a_{0}+a_{1}t+\cdots+a_{k}t^k$. Moreover, the convergence is
uniform on every compact subset of $\mathbb{R}$. Furthermore,
$\lim_{n \to \infty} a_{i,n}=a_i$ for every $i$.
\end{lem}
\begin{proof}[\textbf{Proof of Theorem \ref{main}}] The proof is based on the method of moments. Using the aforementioned properties
of the Kronecker product, the $m^{th}$ moment of the spectral
measure of $\bB_{n,k}$ is given by
\begin{equation}\label{start}
\begin{array}{l c l} \tr_{nk}(\bB_{n,k}^m)&=& \tr_{kn}((\mI_k \otimes \mA_n +
\MW_k \otimes \mB_n)^m) \\
\\
   &=& \sum_{j=0}^m \tr_k(\MW_k^j)\; \tr_n(\phi(\mA_n,\mB_n;m-j,j)),
\end{array}
\end{equation}
where $\phi(\mA_n,\mB_n;m-j,j)$ is the sum of all the noncommutative
monomials in which $\mB_n$ appears $j$ times and $\mA_n$ appears
$m-j$ times.  Let the $j^{th}$ moment of $\mu_{\omega}$ be
$\omega_j$, $j\geq 1$. By \eqref{matmoment} and \eqref{main1},
$\lim_{k \to \infty}\tr_k(\MW_k^{j})=\omega_j$ \emph{a.s.}
Therefore,
\begin{equation}\label{main0}
\lim_{k \to \infty} \tr_{nk}(\bB_{n,k}^m)=\sum_{j=0}^m \omega_j\;
\tr_n(\phi(\mA_n,\mB_n;m-j,j)) .
\end{equation}
On another hand, for all $t\in \mathbb{R}$
\begin{equation}\label{start1}
\tr_n((\mA_n+t\mB_n)^{m})=\sum_{j=0}^m t^j
\tr_n(\phi(\mA_n,\mB_n;m-j,j)).
\end{equation}
Therefore,
\begin{equation*}
\int_{\mathbb{R}} \tr_n((\mA_n+t\mB_n)^{m})
\mu_{\omega}(dt)=\sum_{j=0}^m \omega_j\;
\tr_n(\phi(\mA_n,\mB_n;m-j,j)).
\end{equation*}
So by \eqref{main0},
\begin{equation}\label{main5}
\lim_{k \to \infty} \tr_{nk}(\bB_{n,k}^m)=\int_{\mathbb{R}}
\tr_n((\mA_n+t\mB_n)^{m}) \mu_{\omega}(dt).
\end{equation}
But,
\begin{equation*}
\lim_{n \to \infty}\tr_n((\mA_n+t\mB_n)^{m})=\int_{\mathbb{R}} x^m
\psi(t;dx)
\end{equation*}
and by Lemma \ref{poly} this limit is uniform in $t$ as $t$ varies
over the compact set $\mbox{supp}(\mu_\omega)$. Therefore,
\begin{equation*}
\lim_{n\to \infty} \lim_{k\to
\infty}\tr_{nk}(\bB_{n,k}^m)=\int_{\mathbb{R}} \int_{\mathbb{R}} x^m
\psi(t;dx)\: \mu_{\omega}(dt).
\end{equation*}
By Fubini's Theorem
\begin{equation*}
\lim_{n\to \infty} \lim_{k\to \infty}\tr_{nk}(\bB_{n,k}^m)=
\int_{\mathbb{R}} x^m \left(\int_{\mathbb{R}} \psi(t;dx)\:
\mu_{\omega}(dt)\right).
\end{equation*}
The other iterated limit follows from the observation that
\eqref{start} and \eqref{start1} imply
\begin{equation*}
\tr_{nk}(\bB_{n,k}^m)=\int_{\mathbb{R}} \tr_n((\mA_n+t\mB_n)^m)
\mu_{\MW_k}(dt)
\end{equation*}
for every $n$ and $k$. Since $\mu_{\MW_k}$ is a discrete measure,
\begin{equation}\label{finitek}
\lim_{n\to \infty} \tr_{nk}(\bB_{n,k}^m)=\int_{\mathbb{R}}
\int_{\mathbb{R}} x^m \psi(t;dx) \mu_{\MW_k}(dt).
\end{equation}
By Lemma \ref{poly}, $\int_{\mathbb{R}}  x^m \psi(t;dx)$ is a
polynomial in $t$ and since $\mu_{\MW_k}$ converges in moments, it
follows that
\begin{equation*}
\lim_{k\to \infty} \lim_{n\to \infty}\tr_{nk}(\bB_{n,k}^m)=
\int_{\mathbb{R}} x^m \left(\int_{\mathbb{R}} \psi(t;dx)\:
\mu_{\omega}(dt)\right).
\end{equation*}
Now, since $\nu(dx):=\int_{\mathbb{R}} \psi(t;dx)\:
\mu_{\omega}(dt)$ has a bounded support, the result follows.
\end{proof}
\begin{rem}\label{girko1}
The integral in \eqref{main4} always exists because $\psi(t;.)$ is
measurable in $t$. This can be seen as follows. The characteristic
function of $\psi(t;.)$ is analytic, as $\psi(t;.)$ has a compact
support for each $t$. So  the characteristic function is measurable
in $t$ as a pointwise limit of the series in the moments
$\int_{\mathbb{R}} x^k \psi(t;dx)$,  $k\geq 1$; the latter are
polynomials in $t$ by Lemma \ref{poly}. Therefore the inversion
formula of the characteristic function implies the measurability of
$\psi(t;(-\infty,a])$ for any $a$.
\end{rem}
\section{Applications}
In this section, we apply Theorem \ref{main}  to some well-studied
 ensembles of random matrices.
To do so, we  introduce these ensembles and review the pertinent
topics from free probability theory.
\subsection{Random Matrix Theory}
We call an $n\times n$ Hermitian matrix $\mA=(A_{ij})_{i,j=1}^n$ a
Wigner matrix if it is a random matrix whose upper-diagonal entries
are independent and identically distributed complex random variables
such that $E(A_{ij})=0$ and $E(|A_{ij}|^2)=\frac{1}{n} $ for all
$i<j$. Moreover, the diagonal entries are independent and
identically distributed real random variables such that
$E(A_{ii})=0$ and $E(A_{ii}^2)=\frac{1}{n}$. We will denote all such
Wigner matrices of order $n$ by $\mathbf{Wigner}(n)$.

An important example of a Wigner matrix is the Gaussian Wigner
matrix for which $\{\Re A_{ij}:1 \leq i\leq j\leq n\} \cup \{\Im
A_{ij}:1 \leq i < j\leq n\}$ is a family of independent random
Gaussian variables such that $A_{ii}\sim N(0,\frac{1}{n})$ for every
$i$ and $\Re A_{ij}$, $ \Im A_{ij} \sim N(0,\frac{1}{2n})$ for every
$i < j$. We will denote all such Gaussian matrices of order $n$ by
$\mathbf{G}(n)$.

We call the random matrix $\mB=\mX^* \mX$ a Wishart matrix if $\mX$
is a $p_n \times n$ matrix whose entries are complex independent
Gaussian random variables such that $\Re X_{ij}$, $ \Im X_{ij} \sim
N(0,\frac{1}{2n})$ for every $1 \leq i,j \leq n$. Here $\mX^*$ is
the conjugate transpose of $\mX$. We will denote all such Wishart
matrices of order $n$ and shape parameter $p_n$ by
$\mathbf{Wishart}(n,p_n)$. See \cite{Bai99,Capitaine-Casalis04} for
more details and references.

For these type of random matrices, the limiting spectral
distributions are known. If $\mA_n $ is $\mathbf{Wigner}(n)$, then
by Wigner's Theorem (\emph{cf.} \cite{Bai99}),
$$ \mu_{\mA_n}\xrightarrow{\calD} \gamma_{0,1} \: \mbox{ as
$n\to \infty$} \qquad a.s.$$ where $\gamma_{\alpha,\sigma^2}$ is the
semicircular law centered at $\alpha$ and of variance $\sigma^2$
which is given as
$$\gamma_{\alpha,\sigma^2}(dx)=\frac{1}{2\pi \sigma^2} \sqrt{4 \sigma^2 - (x-\alpha)^2} \:\:
\mathbf{1}_{[\alpha-2\sigma, \alpha+2\sigma]}(x) dx.$$ On the other
hand, if $\mB_n $ is $ \mathbf{Wishart}(n,p_n)$ for all $n$ and
$\lim_{n\to \infty}\frac{p_n}{n}=\alpha>0$, then (\emph{cf.}
\cite{Bai99})
$$\mu_{\mB_n}\xrightarrow{\calD} \rho_\alpha
\: \mbox{ as $n\to \infty$} \qquad a.s.   $$ where $\rho_\alpha$ is
the Marchenko-Pastur law with mean $\alpha>0$ which is given as
\begin{multline*}
\rho_\alpha(dx) =(1-\alpha)^+\delta_0(dx)\\
+\frac{\sqrt{(x-(\sqrt{\alpha}-1)^2)((\sqrt{\alpha}+1)^2-x)}}{2\pi
x}\: \mathbf{1}_{[(\sqrt{\alpha}-1)^2,(\sqrt{\alpha}+1)^2]}(x) dx.
\end{multline*}

The following two propositions are consequences of Theorem
\ref{main}.
\begin{prop}\label{WignerWigner}
Let $\{\MW_n\}$, $\{\mA_n\}$ and $\{\mB_n\}$ be three independent
sequences of random matrices such that $\MW_ n$ is $
\mathbf{Wigner}(n)$ and $\mA_n, \mB_n$ are $ \mathbf{G}(n)$ for all
$n$. Then the almost sure limiting spectral distribution $\nu$ of
$\bB_{n,k}=\mI_k \otimes \mA_n + \MW_k \otimes \mB_n$, see
\eqref{main3}, is absolutely continuous with the probability density
function
\begin{equation*}\label{WWW}
      g(x) = \left \{
        \begin{array}{l@{\quad \: \quad}l}
            g_1(x) & \; ; \text{ whenever } 2 \leq |x| \leq 2 \sqrt{5} \\ \\
             g_2(x)  & \; ; \text{ whenever } |x| \leq 2
           \end{array} \right.
\end{equation*}
where
$$
g_1(x)=\frac{1}{2 \pi^2} \int_{\sqrt{\frac{x^2}{4}-1}}^2
\frac{\sqrt{4(1+t^2)-x^2}\; \sqrt{4-t^2}}{ (1+t^2)} \: dt,
$$
and
$$
g_2(x)=\frac{1}{2 \pi^2} \int_0^2 \frac{\sqrt{4(1+t^2)-x^2}\;
\sqrt{4-t^2}}{ (1+t^2)} \:dt.
$$
\end{prop}
\begin{figure}[hbt]
\includegraphics[height=3.5cm]{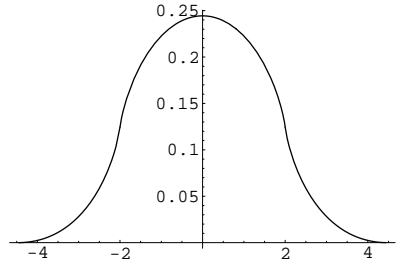}
\caption{The probability density function corresponding to the
limiting spectral distribution of $\bB_{n,k}$ when $\mA_n$ and
$\mB_n$ are Gaussian matrices. }
\end{figure}

In order to state the following proposition we first define the
following functions:
$$h_1(x;t)=2+27\,t^2-3\,t\,x - 3\,t^2\,x^2+ 2\,t^3\,x^3 ,$$
$$h_2(x;t)=1 - t\,x + t^2\,x^2,$$
$$H(x;t)=\frac{1}{\sqrt[3]{2}} \left( h_1(x;t) +\sqrt{h_1(x;t)^2-4\,
h_2(x;t)^3 }\right) ^{\frac{1}{3}},$$  the two functions $s_1(t)$
and $s_2(t)$ which are the two real roots of the quartic equation in
$x$
\begin{equation}\label{intineq}
4 + 27\,t^2 - 6\,t\,x - x^2 - 6\,t^2\,x^2 + 2\,t\,x^3 + 4\,t^3\,x^3
- t^2\,x^4=0
\end{equation}
(see Proposition \ref{aux2} for details) and the probability density
function
\begin{equation}\label{a+tb}
      f(x;t) = \frac{1}{2\,\sqrt{3}\,\pi \,t}\left(
H(x;t)-\frac{h_2(x;t)}{H(x;t)} \right) \;
\mathbf{1}_{[s_1(t),s_2(t)]}(x).
\end{equation}
\begin{prop}\label{WignerWishart}
Let $\{\MW_n\}$, $\{\mA_n\}$ and $\{\mB_n\}$ be three independent
sequences of random matrices such that $\MW_n $ is
$\mathbf{Wigner}(n)$, $\mA_n $ is $ \mathbf{G}(n)$ and $\mB_n$ is $
\mathbf{Wishart}(n,p_n)$ for all $n$ and $\lim_{n\to
\infty}\frac{p_n}{n}=1$. Then the almost sure limiting spectral
distribution  $\nu$ of $\bB_{n,k}=\mI_k \otimes \mA_n + \MW_k
\otimes \mB_n$, see \eqref{main3}, is absolutely continuous with the
probability density function
\begin{equation*}
g(x)=\int_{\mathbb{R}} f(x;t)\: \gamma_{0,1}(dt)
\end{equation*}
where $f(x;t)$ is given in Equation \eqref{a+tb}.
\end{prop}

We simply need to verify the hypothesis of Theorem \ref{main}. We
use the machinery of free probability (specifically free additive
convolution) to do this.
\subsection{Free Additive Convolution}
Let $\mu$ be a probability measure with a compact support in
$\mathbb{R}$. We define its corresponding Cauchy (Stieltjes)
transform to be $G_{\mu}(z):=\int_{\mathbb{R}}\frac{1}{z-x}\mu(dx)$,
for $z\in \mathbb{C}$ such that $\Im (z) >0$. The Cauchy transform
$G_{\mu}(z)$ possesses the following properties:
\begin{enumerate}
\item  $\Im (G_{\mu}(z))<0$ whenever $\Im (z)>0$.
\item $\lim_{|z|\to \infty} z G_{\mu}(z)=1$.
\item $G_{\mu}(z)$ is analytic in a neighborhood of $\infty$.
\end{enumerate}
The Stieltjes inversion formula is given by
\begin{equation}\label{inv}
\mu(E)=- \frac{1}{\pi} \lim_{y\downarrow 0} \int_{E}
\Im(G_{\mu}(x+\mathbf{i} \,y)) dx
\end{equation}
for every continuity set $E\subset \mathcal{B}(\mathbb{R})$ (the
$\sigma-$field of Borel subsets of $\mathbb{R}$). The R-transform of
$\mu$ is defined as $R_{\mu}(z)=K_{\mu}(z)-\frac{1}{z}$ where
$K_{\mu}(z)$ is the inverse function of the Cauchy transform
$G_{\mu}(z)$, \emph{i.e.}, $G_{\mu}(K_{\mu}(z))=z$. The two
functions $K_{\mu}(z)$ and $R_{\mu}(z)$ are well defined in
$0<|z|<r$ and $0 \leq |z|<r$, respectively, for some $r>0$.

The free additive convolution of probability measures with compact
supports in $\mathbb{R}$ arises in free probability theory
(\emph{cf.} \cite{VDN92}). If $\mu$ and $\nu$ are two probability
measures with compact supports in $\mathbb{R}$, then their free
additive convolution $\mu \boxplus \nu$ is a probability measure
with a compact support in $\mathbb{R}$, see \cite{VDN92}. The
R-transform of $\mu \boxplus \nu$ is given by $R_{\mu \boxplus
\nu}(z)=R_{\mu}(z)+R_{\nu}(z)$.

Denote the dilation $D_t$ of a measure $\mu$ by $D_t(\mu)$, where
$D_t(\mu)(E)=\mu(E/t)$ for every $E$ if $t \neq 0$ and
$D_t(\mu)=\delta_0$ if $t = 0$. Since $R_{D_t(\mu)}(z)=t R_{\mu}(t
z)$ for every $t \in \mathbb{R}$ (\emph{cf.} \cite[p.26]{VDN92}),
therefore
$$D_t(\mu \boxplus \nu)=D_t(\mu)\boxplus D_t(\nu)
.$$ For the limit laws described above,
$R_{\gamma_{0,\sigma^2}}(z)=\sigma^2 z$ and
$R_{\rho_\alpha}(z)=\frac{\alpha}{1-z}$ (\emph{cf.} \cite{VDN92});
furthermore $\gamma_{0,1} \boxplus
D_t(\gamma_{0,1})=\gamma_{0,1+t^2}$ for every real $t$. The next
proposition computes $\gamma_{0,1}\boxplus D_t(\rho_1)$ for $t>0$.
\begin{prop}\label{aux2}
If $\mu=\gamma_{0,1}$ and $\nu=\rho_1$ then for $t>0$ $$\mu \boxplus
D_t(\nu) (dx)=f(x;t) dx$$ where $f(x;t)$ is given by Equation
\eqref{a+tb}. Furthermore, for each compact set $C \subset
\mathbb{R}$ there exists $M>0$ such that the support of $f(x;t) dx$
is contained in $[-M,M]$ for all $t$ in $C$.
\end{prop}
\begin{proof}
Fix $t >0$. The R-transforms of $\mu$ and $\nu$ are given by
$$R_{\mu}(z)= z \text{ and
}R_{D_t(\nu)}(z)=\frac{t}{1-tz}$$ and accordingly
$$R_{\mu \boxplus D_t(\nu)}(z)= z+\frac{t}{1-tz}.$$ Therefore, the Cauchy
transform $G_{\mu \boxplus D_t(\nu)}(z)$ is the root of the cubic
equation
\begin{equation}\label{cubic}
t\,g^3 - g^2\,\left( 1 + t\,z \right)+ g\,z -1 = 0.
\end{equation}
First, in order to show uniqueness, we will show that Equation
\eqref{cubic} has only one root for which $\lim_{|z|\to \infty} z
g(z)=1$. This follows from the observation that if $g_1$, $g_2$ and
$g_3$ are the roots of Equation (\ref{cubic}) then
$g_1g_2g_3=\frac1t$ and $g_1+g_2+g_3=\frac1t+z$. Combining both
identities results in
\begin{equation}\label{combine}
g_1^2 g_2+g_1g_2^2+\frac1t=\frac1t g_1 g_2+zg_1 g_2
\end{equation}
Thus if two of the roots, say $g_1$ and $g_2$, are such that
$\lim_{|z|\to \infty} z g_1(z)=1$ and $\lim_{|z|\to \infty} z
g_2(z)=1$, then Equation \eqref{combine} would lead to the
contradiction that $\frac1t=0$. It is known, see \cite[Corollary 2,
Corollary 4, and Proposition 5]{Biane97}, that the free convolution
of a compactly supported measure with a semicircle law has a smooth
and bounded density. Thus by picking the right root and then using
the Stieltjes inversion formula \eqref{inv} we get
\begin{equation*}
- \frac{1}{\pi} \lim_{y\downarrow 0} \Im(G_{\mu \boxplus
D_t(\nu)}(x+\mathbf{i} \,y)) =f(x;t)
\end{equation*}
where $f(x;t)$ is given by Equation (\ref{a+tb}).

Second, since we know in advance that the free additive convolution
of two probability measures with compact supports in $\mathbb{R}$
has a compact support in $\mathbb{R}$ (\emph{cf.} \cite{VDN92}),
then we can find the support of $\mu \boxplus D_t(\nu) $ by
identifying when $f(x;t)=0$. This last identity leads to Equation
(\ref{intineq}).

The left hand side in Equation (\ref{intineq}) is positive at $x=0$
and Equation (\ref{intineq}) has two real roots and two complex
conjugate roots. We prove the existence and uniqueness of its real
roots $s_1(t)$ and $s_2(t)$ as follows. Substitute $x=y+\frac{1}{2t}
+ t$ in Equation (\ref{intineq}). Hence, we get another quartic
equation in $y$ which reads as
\begin{equation}\label{quartic2}
- \left(\frac{37}{2} -
  \frac{1}{16\,t^4} + 3\,t^4 \right)+
 \left( \frac{8}{t} - 8\,t^3 \right)\,y -
  \left(\frac{1}{2\,t^2} + 6\,t^2\right)\,y^2 +
  y^4=0
\end{equation}

By Descartes$-$Euler theorem \cite{Abra}, Equation \eqref{quartic2}
has two real roots and two complex conjugate roots (and
correspondingly Equation (\ref{intineq}) does) if and only if the
cubic equation in $z$ that is given by
\begin{equation}\label{cubic2}
-{\left( \frac{8}{t} - 8\,t^3 \right)
       }^2 + \left( 80+48 \,t^4
        \right)   \,z -
  \left( \frac{1}{t^2} +
     12\,t^2 \right) \,z^2 + z^3=0
\end{equation}
has one nonnegative real root and two complex conjugate roots. To
show this we substitute $z=u+\left( \frac{1}{3\,t^2} +
       4t^2 \right) $ in Equation \eqref{cubic2} which in its turn
reduces to another cubic equation. By Cardano's theorem, the
discriminant of the new reduced equation $D=\frac{64\,{\left( 1 +
27\,t^4 \right)
        }^3}{27\,t^8}$ is nonnegative and so it has one real root and two complex conjugate
roots and so does Equation (\ref{cubic2}).

Now, it is left to show that this real root of Equation
(\ref{cubic2}) is nonnegative. This is true since the product of the
three roots of Equation (\ref{cubic2}) is equal to ${\left(
\frac{8}{t} - 8\,t^3 \right)}^2$ which is nonnegative and
consequently the real root is also nonnegative.

We can see from Equation (\ref{intineq}) that $s_1(t)$ and $s_2(t)$
are continuous and hence uniformly bounded on any compact subset of
$\mathbb{R}$. This completes the proof.
\end{proof}
In Figure 1, we show the graphs of $f(x;t)$ for $t=1,2$.
\begin{figure}[hbt]
\begin{tabular}{ccc}
\includegraphics[height=3.5cm]{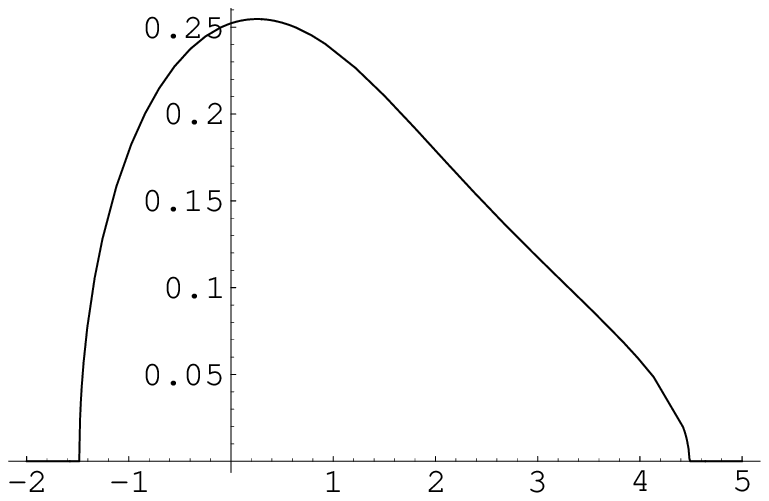}& \includegraphics[height=3.5cm]{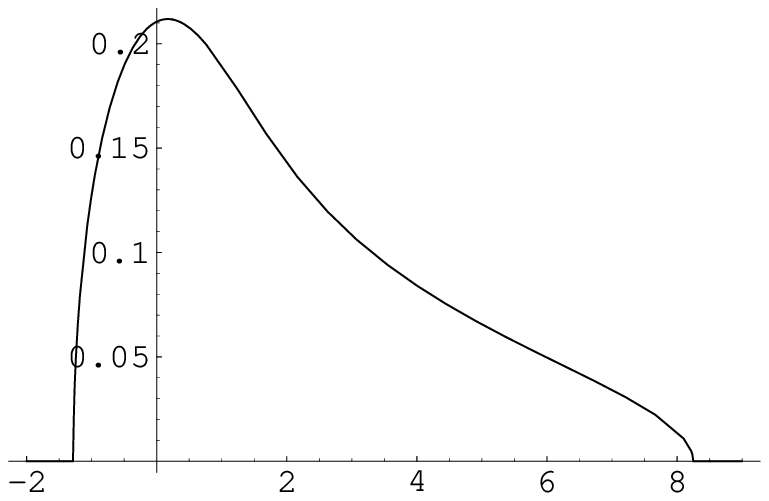}\\
$f(x;1)$ & $f(x;2)$  \\
\end{tabular}
\caption{The probability density functions corresponding to \newline
$\mu \boxplus D_t(\nu)$ with t=1, 2. \label{F1}}
\end{figure}
\begin{rem}\label{tzero}
If $t=0$, then $f(x;0)=\frac{1}{2 \pi} \sqrt{4-x^2}\;
\mathbf{1}_{[-2,2]}(x)$ since $\mu \boxplus \delta_0=\mu$. In case
$t<0$, since $\mu \boxplus D_t(\nu)=D_{-1}(\mu \boxplus
D_{-t}(\nu))$, then it follows that $\mu \boxplus
D_t(\nu)(dx)=f(-x;-t)dx$.
\end{rem}
We will also need the following.
\begin{thm}\label{freecon}\cite[Proposition 4.3.9]{Hiai-Petz}
Let $\{\mA_n\}$ and $\{\mB_n\}$ be two sequences of Hermitian random
matrices and $\{\mU_n\}$ be a sequence of random matrices with the
uniform distribution on the unitary group $\mathbb{U}(n)$. Suppose
that $\mU_n$ is independent of $(\mA_n,\mB_n)$ for all $n\geq 1$. If
there exist two compactly supported probability distributions $\mu$
and $\nu$ such that
$$\mu_{\mA_n}\xrightarrow{\frakm}
\mu \: \mbox{   as $n\to \infty$} \qquad a.s.\: \text{ and } \:
\mu_{\mB_n}\xrightarrow{\frakm} \nu \: \mbox{   as $n\to \infty$}
\qquad a.s., $$ then
$$\mu_{\mA_n+\mU_n^*\mB_n\mU_n}\xrightarrow{\calD}
\mu\boxplus\nu \: \mbox{   as $n\to \infty$} \qquad a.s.$$
\end{thm}

In particular, if $\mA_n$, $\mB_n$ are independent, and the
distribution of $\mB_n$ is unitarily invariant for all $n$,
\emph{i.e.}, if $\mB_n$ and $\mU_n \mB_n \mU_n^*$ have the same
distribution for all unitary matrices $\mU_n$, then Theorem
\ref{freecon} implies
$$\mu_{\mA_n+\mB_n}\xrightarrow{\calD}
\mu\boxplus\nu \: \mbox{   as $n\to \infty$} \qquad a.s.$$ Since the
distributions of $\mathbf{G}(n)$ and $\mathbf{Wishart}(n,p_n)$
matrices are unitarily invariant, we get the following.
\begin{cor}\label{A+tB}
\noindent
\begin{enumerate}
\item For $n\geq 1$, let $\mA_n$ be $ \mathbf{Wigner}(n)$ and
$\mB_n $ be $ \mathbf{G}(n)$ such that $\mA_n$ and $\mB_n$ are
independent. Then for all $t \in \mathbb{R}$,
\begin{equation*}
\mu_{\mA_n+t\mB_n}\xrightarrow{\calD} \gamma_{0,1+t^2} \: \mbox{ as
$n\to \infty$} \qquad a.s.
\end{equation*}
\item For $n\geq 1$, let $\mA_n $ be $ \mathbf{Wigner}(n)$ and
$\mB_n $ be $ \mathbf{Wishart}(n,p_n)$ such that $\mA_n$ and $\mB_n$
are independent. If $\lim_{n\to \infty}\frac{p_n}{n}=1$,
 then for all $t \in \mathbb{R}$,
\begin{equation*}
\mu_{\mA_n+t\mB_n}\xrightarrow{\calD} \mu_t \: \mbox{   as $n\to
\infty$} \qquad a.s.
\end{equation*}
where $\mu_t$ has the probability density function given by Equation
(\ref{a+tb}).
\end{enumerate}
\end{cor}

Now, Proposition \ref{WignerWigner} follows directly from Corollary
\ref{A+tB} part (i). Proposition \ref{WignerWishart} follows easily
from Corollary \ref{A+tB} part (ii) since $f(x;t)dx$ has a bounded
support that is uniformly bounded in $t \in [-2,2]$ as shown in
Proposition \ref{aux2}.

\subsection{Comments on limits with respect to one index} Here we remark about additional information about limits with respect to one of the indexes
that can be extracted from the proof of Theorem \ref{main}.
\subsubsection{}
Formula \eqref{finitek} identifies the limiting spectral
distribution of finite dimensional random block-matrices as the size
of the blocks goes to infinity. For instance, consider the sequence
of $k \times k$ random block-matrices $\{\bS_n\}$, for a fixed $k
\in \mathbb{N}$, in which the diagonal blocks are made of $\mA_n$'s
and all the other blocks are made of $\mB_n$'s. This random block
matrix is studied in \cite{Tamer2005} among other things, using
algebraic manipulations of block matrices. Hence, for each $n$ we
can write $\bS_n$ as $\mI \otimes \mA_n + \MW \otimes \mB_n$, where
$\mI$ is the $k \times k$ identity matrix and $\MW$ is the $k \times
k$ non-random matrix whose entries are 0's on the diagonal and one's
elsewhere. By induction on $k$, one can easily find that
$\mu_{\MW}=\frac{k-1}{k} \delta_{-1}+\frac1k \delta_{k-1}$. If
$\{\mA_n\}$ and $\{\mB_n\}$ satisfy condition \eqref{2)} in Theorem
\ref{main}, then it follows from formula \eqref{finitek} that
\begin{equation*}
\mu_{\bS_n} \xrightarrow{\calD} \frac{k-1}{k}\psi(-1,.)+\frac{1}{k}
\psi(k-1,.) \: \mbox{ as $n\to \infty$} \qquad a.s.
\end{equation*}

\subsubsection{}
The following result is a generalization of \cite[Theorem
4]{Girko2000}, where it is assumed that the entries of $\MW_k$ take
values $-1$ and $1$ with probabilities $\frac12$, and $\mA$ and
$\mB$ are non-random matrices.
\begin{prop}\label{commute}{\bf{(The SS-Law)}}
Let $\{\MW_k\}$ be a sequence of random matrices such that $\MW_k $
is $ \mathbf{Wigner}(k)$ for all $k$. Let also $\mA$ and $\mB$ be
two $n \times n$ Hermitian commuting random matrices which have
eigenvalues $\{\alpha_1 , \cdots , \alpha_n\}$ and $\{\beta_1 ,
\cdots , \beta_n\}$, respectively. Then
\begin{equation}\label{last}
\mu_{\bB_{n,k}}\xrightarrow{\calD} \frac1n \sum_{j=1}^n
\gamma_{\alpha_j,\beta_j^2} \: \mbox{ as $k\to \infty$} \qquad a.s.
\end{equation}
\end{prop}

\begin{proof} Since $\mA$ and $\mB$ are Hermitian and commute, then
there is a unitary matrix $\mU$ such that
$$\mA=\mU\,\mbox{\bf{Diag}}(\alpha_1,\ldots,\alpha_n) \,\mU^* \text{
and } \mB=\mU\,\mbox{\bf{Diag}}(\beta_1,\ldots,\beta_n)\, \mU^*.$$
It follows directly from equation \eqref{main5} that
\begin{equation*}
\begin{array}{l c l} \lim_{k\to \infty}\tr_{nk}(\bB_{n,k}^m)&=&\int_{\mathbb{R}} \frac1n
\sum_{j=1}^n(\alpha_j+t\beta_j)^m \gamma_{0,1}(dt)  \\
\\
   &=& \frac1n \sum_{j=1}^n \int_{\mathbb{R}} t^m  \gamma_{\alpha_j,\beta_j^2}(dt)
\end{array}
\end{equation*}
where the last equality follows by changing of variables.

Since the probability measure on the right-hand side of \eqref{last}
is a finite mixture of semicircle laws, then it has a compact
support in $\mathbb{R}$. Therefore, convergence of moments implies
weak convergence.
\end{proof}

\section*{Acknowledgement}
I am indebted to my advisor Prof. W. Bryc for his continuous and
sincere help and support. I would like to thank Prof. C. David Minda
for the useful discussions I had with him. I am also thankful for
Prof. J. Mitro for the useful comments. I would like also to thank
the referee and the associate editor for their helpful comments and
suggestions that brought the paper to its final shape.
\bibliographystyle{plain}

\end{document}